\documentclass[10pt]{article}
\usepackage{amsfonts,amsmath,amsthm,latexsym,amssymb,marvosym,esint,xfrac,bbm,amssymb}
\usepackage{mathrsfs}
\usepackage{graphics, epsfig}
\usepackage{color}
\usepackage{appendix}
\usepackage{verbatim}
\usepackage{ulem}
\usepackage[makeroom]{cancel}
\usepackage[margin=1in]{geometry}
\usepackage{tikz}
\tikzstyle{mybox} = [draw=black, very thick, rectangle, rounded corners, inner ysep=5pt, inner xsep=5pt]
\usepackage{pgfplots}

 \usepackage[usenames,dvipsnames]{pstricks}
 \usepackage{pst-grad} % For gradients
 \usepackage{pst-plot} % For axes
\allowdisplaybreaks

\usepackage[colorlinks=true,linkcolor=blue,citecolor=red]{hyperref}% for hyperlinks of references

\newtheorem{proposition}{Proposition}[section]
\newtheorem{theorem}{Theorem}[section]
\newtheorem{definition}{Definition}[section]
\newtheorem{lemma}{Lemma}[section]

\newtheorem{remark}{Remark}[section]

\numberwithin{equation}{section}
\numberwithin{theorem}{section}
\numberwithin{proposition}{section}
\numberwithin{lemma}{section}
\numberwithin{remark}{section}
\newcommand{\essup}{\operatornamewithlimits{ess\,sup}}
\newcommand{\essinf}{\operatornamewithlimits{ess\,inf}}

\newcommand{\intl}{\int\limits}

\def\Xint#1{\mathchoice
    {\XXint\displaystyle\textstyle{#1}}%
    {\XXint\textstyle\scriptstyle{#1}}%
    {\XXint\scriptstyle\scriptscriptstyle{#1}}%
    {\XXint\scriptscriptstyle\scriptscriptstyle{#1}}%
    \!\int}
\def\XXint#1#2#3{\setbox0=\hbox{$#1{#2#3}{\int}$}
    \vcenter{\hbox{$#2#3$}}\kern-0.5\wd0}
\def\bint{\Xint-}
\def\dashint{\Xint{\raise4pt\hbox to7pt{\hrulefill}}}
\def\dashiint{\bint\kern-0.15cm\bint}
\newcommand{\ovl}[3]{\int_{#1}^{#2}\kern-#3pt\raise4pt\hbox to7pt{\hrulefill}\ }

\newcommand{\ovll}[3]{\intl_{#1}^{#2}\kern-#3pt\raise4pt\hbox to7pt{\hrulefill}\ }

\newcommand{\tvl}[2]{\iint_{#1}\kern-#2pt\raise4pt\hbox to7pt{\hrulefill}\ }

\newcommand{\bye}{\end{document}}

\newcommand{\tvls}[2]{\iint_{#1}\kern-#2pt\raise4pt\hbox to15pt{\hrulefill}\ }

\def\R{\mathbb{R}}
\def\Q{\mathcal{Q}}
\def\K{\mathcal{K}}
\def\N{\mathbb{N}}
\def\dive{\mathrm{div}}

 \def\B{\mathcal{B}}

\def\L{\mathcal{L}}
\def\E{\mathcal{E}}

\def\T{\mathcal{T}}
\def\b{\mathbb{B}}
\def\l{\mathscr{L}}
\date{}
\begin{document}
\title{On a particular scaling for the prototype anisotropic p-Laplacian}
%%%%%%%%%%%%%%%%%%%%%%%%%%%%%%%%%
%%%%%%%%%%%%%%%%%%%%%%%%%%%%%%%%%%%%%%%%%%%%%%%%%%%
\author{
\\
\\
\it{Simone Ciani \& Umberto Guarnotta \&  Vincenzo} Vespri \\ \\
Technische Universit\"at Darmstadt, Department of Mathematics, \\
Schlossgartenstraße 7,
64289, Darmstadt,
Germany\\
{\it ciani@mathematik.tu-darmstadt.de}\\ \\
Universit\`a degli Studi di Palermo, Dipartimento di Matematica e Informatica,\\ Via Archirafi 34, 90123, Palermo, Italy\\
{ \it umberto.guarnotta@unipa.it} \\ \\
Universit\`a degli Studi di Firenze, Dipartimento DIMAI\\ Viale G. Morgagni 67/a, 50134, Firenze, Italy \\  { \it vincenzo.vespri@unifi.it}
}
%%%%%%%%%%%%%%%%%%%%%%%%%%%%%%%%%%%%%%%%%%%%%%%%%%%
%%%%%%%%%%%%%%%%%%%%%%%%%%%%%%%%%%%%
\maketitle 
\vskip.4truecm 
%%%%%%%%%%%%%%%%%%%%%%%%%%%%%%%%%%%%%%%%%%%%%%%%%%%
%%%%%%%%%%%%%%%%%%%%%%%%%%%%%%%%%%%%%%%%%%%%
\begin{abstract} \noindent
%%%%%%%%%%%%%%%%%%%%%%%%%%%%%%%%%%%%%%%%%%%%%%%%%%%
\vskip.2truecm
%%%%%%%%%%%%%%%%%%%%%%%%%%%%%%%%%%%%%%%%%%%%%%%%%%% 
\noindent 
In this brief note we show that under a volume non-preserving scaling it is possible to recover the basics for a regularity theory regarding local weak solutions to the fully anisotropic equation 
\begin{equation} \label{EQ}
    \partial_t u= \sum_{i=1}^N \partial_i (|\partial_i u|^{p_i-2} \partial_i u) \quad \text{in} \quad \Omega_T= \Omega \times (-T,T), \quad \text{with} \quad \Omega \subset \subset \R^N.
\end{equation} \noindent We characterize self-similar solutions regarding this particular scaling and we show that semi-continuity for solutions to this equation is a consequence of a simple property that is itself invariant under scaling.
\vskip.2truecm
%%%%%%%%%%%%%%%%%%%%%%%%%%%%%%%%%%%%%%%%%%%%%%%%%%%
\noindent
{\bf{MSC 2020:}}
35K65, 35K92, 35B65.\vskip0.5cm
\noindent
{\bf{Keywords}}: Anisotropic $p$-Laplacian, Critical Mass Lemma, Intrinsic Scaling, Lower Semi-Continuity. 
%%%%%%%%%%%%%%%%%%%%%%%%%%%%%%%%%%%%%%%%%%%%%%%%%%%

\begin{flushright}
\it{To celebrate Francesco Altomare's 70th genethliac}
\end{flushright}
\end{abstract}
%%%%%%%%%%%%%%%%%%%%%%%%%%%%%%%%%%%%%%%%%%%%%%%%%%%%%%%%%%%%%%%%%%%%%%%

\section{Introduction to the problem}
Equation \eqref{EQ} is a parabolic anisotropic equation with non-standard growth. We refer to the introduction of \cite{Brasco1}, \cite{UU} and the surveys \cite{Marcellini}, \cite{Mingione}  for a non-exhaustive introduction to the origin of the problem, and to the introduction of \cite{Brasco2} and the book \cite{Antontsev-Shmarev} for a more general account to the parabolic problem. At a first glance equation \eqref{EQ} may look similar to the equation \begin{equation}\label{pLaplace}
u_t- \dive (|\nabla u|^{p-2} \nabla u)=0 \quad \text{locally weakly in}\, \, \Omega_T.
\end{equation}

\noindent Literature on this topic is very developed, and even if the problem of regularity of solutions to \eqref{EQ} is old more than fifty years, still very much is unknown from the point of view of basic regularity, as local H\"older continuity or Harnack inequality. The principal motivation is that the techniques usually employed for nonlinear equations (as $p$-Laplacian equations, porous medium equations, doubly nonlinear equations, and so on) are not directly applicable to it. Let us explain this point in detail.\vskip0.2cm \noindent
Up to our knowledge, in the setting of evolutionary nonlinear operators of $p$-growth (whose prototype is \eqref{pLaplace} with $p\neq2$), the main technique to prove a Harnack inequality is exploiting a parabolic continuous transformation having the general form 
\begin{equation} \label{giorgino} w(x,t)= e^{t/(p-2)} u(x,e^t), \quad \quad  x \in \Omega, \, t>0. \end{equation}\noindent This transformation maps super-solutions to \eqref{pLaplace} to super-solutions to a similar equation, that has an exponential dependence on time only on the non-homogeneous terms. Along this strategy, the possibility to stretch time and control the non-homogeneous terms is crucial, in order to employ a technique originally conceived by E. DeGiorgi for solutions to elliptic partial differential equations (see, e.g., \cite{DG}, \cite{acta}), based on a version of the isoperimetric inequality (cf. \cite[Lemma 2.2., page 5]{DB}). This argument allows to prove an expansion of positivity for the transformed super-solutions that, if carried back to solutions to \eqref{pLaplace}, provides the expansion of positivity necessary for an intrinsic Harnack inequality to hold true.\vskip0.2cm \noindent 
The main issue dealing with \eqref{EQ} is that, in general, a continuous transformation with an exponential-type dependence on time necessarily affects the space variables. Taking into account also the strong nonlinear behavior of the equation along the space variables, the control of the non-homogeneous terms in the transformed equation is encumbered. More precisely, from the energetic point of view, the new equation is no more of the same kind of \eqref{EQ}, and this leads the whole machinery to fail.\vskip0.2cm \noindent On the other hand, in \cite{Ciani-Mosconi-Vespri} the authors proved that an intrinsic Harnack type inequality is valid for local weak solutions to \eqref{EQ}, by adapting a classic idea of E. DiBenedetto (see \cite{DB2}) consisting in a comparison between the solution and a particular one, called Barenblatt solution in honor to its discoverer (see the original in \cite{Barenblatt} and \cite{Ciani-Vespri} for an overview on the anisotropic case). However, the generalization of this inequality to a wider class of parabolic operators patterned after \eqref{EQ} is still an open and challenging problem. The purpose of the present work is to investigate a particular scaling of the equation: it would permit to free the time variable from the space ones, opening the way to an application of a transformation similar to \eqref{giorgino}. This homogenization seems to unveil a new insight on the anisotropic behavior of these operators. From the energetic point of view, serious difficulties appear even with the stationary counterpart of \eqref{EQ}, because the competition among different directional $p_i$-diffusions encodes both singular and degenerate behavior. Roughly speaking, this can be illustrated within the scaling of \cite{Ciani-Skrypnik-Vespri}, looking at the kind of degeneration that the set 
\[\prod_{i=1}^N \bigg{\{} |x_i|< \rho^{\frac{\bar{p}}{p_i}} M^{\frac{p_i-\bar{p}}{p_i}}  \bigg{\}}, \qquad M,\rho>0\]
exhibits as $M$ vanishes. This is a volume-preserving set of self-similar geometry where the equation evolves, and the parameter $M$ is usually chosen to be a multiple of the oscillation of $u$, in order to restore the homogeneity of the energy. The problem is that, depending on the sign of $(p_i-\bar{p})$, the set stretches or vanishes along the respective coordinates.  \vskip0.2cm 

% \noindent We propose here a different intrinsic scaling of the equation (see \eqref{Trans}-\eqref{transf}), that enjoys the following properties. First of all, the time time variable is free from intrinsic geometry; this means that it is not linked to the solution itself. Secondly, this geometry associated with the scaling degenerates monotonically with $M$, so we say that the geometry is only {\it degenerate}, not singular. This means that as the solution gets smaller then the intrinsic geometry shrinks in every direction according to it.

% \vskip0.4cm 
\noindent The different scaling that we propose in this note (see \eqref{Trans}-\eqref{transf}) possesses the following properties: the intrinsic geometry associated with it degenerates monotonically with $M$, so we say that the geometry is only {\it degenerate}, not singular; it does not affect anyhow the time variables from the intrinsic point of view. From this perspective, this particular scaling seems promising; see for instance, the energy $\E_n$ in Lemma \ref{DG}. The crucial point is that we can identify the self-similar solutions to \eqref{EQ}, namely, the solutions that coincide with their scaled functions; this is done via correspondence with a Fokker-Planck equation (cf. Proposition \ref{fk-ss}). As a consequence, all the properties of solutions to \eqref{EQ} proved in \cite{Ciani-Mosconi-Vespri} hold true, in a re-interpreted formulation, also for solutions of a `wild' Fokker-Planck equation (see \eqref{FokkerPlanck}). The existence of a Barenblatt solution is of fundamental importance to understand the behaviour of solutions.

\vskip0.2cm \noindent 
Moreover, we show that this special scaling preserves the energy of the solutions, as well as other properties, that will be called for this reason {\it invariants} (standing for scale-invariants). An example is furnished by the Critical Mass Lemma, that can be regarded as a measure-theoretical maximum principle (see Lemma \ref{DG} for details; see also \cite[p. 8]{VespriINDAM} and Proposition 2 of \cite{Eurica} in the context of anisotropic porous medium).
\vskip0.2cm \noindent 
Dimensional analysis is a simple consequence of the well-known covariance principle of Physics: all physical laws can be represented in a form which is equally valid for all observers. The very idea of self-similarity is connected with the group of transformations of solutions: see, e.g., \cite{Barenblatt-scaling}. These groups are inborn in the differential equations governing the process, and are
determined by the physical dimensions of the variables appearing in them; transformations of units of time, length, mass, etc. are the simplest examples. This kind of self-similarity is obtained by power laws with exponents that are simple fractions defined in an elementary way from dimensional considerations. These arguments led to an interpretation of nonlinear parabolic theory, developed by DiBenedetto (\cite{DB}), Vazquez (\cite{Vazquez}), and many others, which is nowadays known as method of intrinsic scaling (cf. also \cite{Urb}). The key feature of the argument of intrinsic scaling is that, by appropriately scaling the geometry, the energy of solutions enjoys a homogeneous form that is easier to manipulate. This idea can be used in turn to interpret the energy of solutions to anisotropic equations like \eqref{EQ} in a homogeneous fashion. This is the purpose of the present scaling, whose side-effect on energy $\E_n$ is here shown by the non-scaled version of Lemma \ref{De Giorgi}. We present here a general version of this lemma, that we could not find in literature for the full parabolic anisotropic equation and that is propaedeutic to the study of further properties.

\vskip0.2cm 

\noindent Indeed, as a byproduct of our analysis, by applying the ideas of \cite{Liao} to the parabolic setting, we show that lower semi-continuity of super-solutions is a sole consequence of these general invariants. \vskip0.2cm \noindent 
The existence of a lower semi-continuous representative for local weak super-solutions has already been obtained in \cite{Mosconi} by using an idea of \cite{Kuusi}. The authors observe that a proper $L^r-L^{\infty}$ estimate for weak super-solutions suffices to obtain a lower semi-continuous representative. This technique is however linked to the particular structure of the equation, that allowed them to add a constant to the solution to generate another solution. The new approach of \cite{Liao} is more general, since the existence of a lower semi-continuous representative is linked only to a more general property, that is the analogue of Lemma \ref{DG}. In this way the authors of the aforementioned \cite{Mosconi}, \cite{Kuusi}, \cite{Liao} proved that weak solutions are $p$-super-harmonic solutions. The latter ones are, on an appropriate setting, proper lower semi-continuous functions, that can be compared with any sufficiently regular solution to the same equation. Since the comparison principle for equations driven by monotone operators holds true, the main step consists in proving semi-continuity.\newline 
It would be an interesting subject to determine whether $p$-super-harmonic functions, whose derivatives a priori may be even unbounded, can satisfy a Critical Mass Lemma as Lemma \ref{DG}. 

\subsection*{Structure of the paper}
Section \ref{Preliminaries} is devoted to set up the functional framework and propose the particular scaling. In Section \ref{invariants} we show that energetic properties of the equation are invariant under this scaling. Finally, in Section \ref{Ssemicontinuity}, we furnish a new proof of semi-continuity for super-solutions to \eqref{EQ}.

\subsection*{Notation}
{\small
\begin{itemize}
    \item[-] Let $N\in\N$, $N>1$, and let $\Omega\subset \R^N$ be a bounded open set. Given $T>0$, we set $\Omega_T=\Omega \times (-T,T)$. \\
    The symbol $A \subset\subset B$ means `$A$ is compactly contained in $B$'.
    
    \item[-] For any $\varphi \in W^{1,1}(\Omega)$, we denote by $\partial_i \varphi$ the $i$-th directional weak derivative of $\varphi$. \\
    If moreover $\varphi \in W^{1,2}([s,t], L^2(\Omega))$ for some $s,t\in\R$, $s<t$, then $\partial_t \varphi$ stands for the weak time-derivative of $\varphi$.
    
    \item[-] We denote the cube of side $2\rho>0$ and center $x \in \R^N$ with $x+K_{\rho}$, while \newline
    \begin{equation*} \begin{cases} (x,t)+ Q_{\rho,\tau}^-= (x+K_{\rho}) \times (-\tau,0],\\
    (x,t)+ Q_{\rho,\tau}^+= (x+K_{\rho}) \times [0,\tau),\\
    (x,t)+ Q_{\rho,\tau}= (x+K_{\rho}) \times (-\tau,\tau], \end{cases} \end{equation*} stand for, respectively, the backward, forward and full cylinders centered at $(x,t) \in \R^{N+1}$. \\
    When $\tau=1$ we simply write $Q_\rho^-$,$Q_\rho$,$Q_\rho^+$ instead of $Q_{\rho,1}^-,Q_{\rho,1},Q_{\rho,1}^+$.
    
    \item[-] We fix a vector of $N$ numbers ${\bf{p}}=(p_1,\dots,p_N)$; the index $i$ will run through $1,\ldots,N$. We define the harmonic mean of $p_i$s as $\bar{p}= N(\sum_{i=1}^N 1/p_i)^{-1}$, and for $\bar{p}<N$ the Sobolev exponent of the harmonic mean by $\bar{p}^*= N\bar{p}/(N-\bar{p})$. Hereafter we suppose $$2<p_1\leq p_2\leq \dots \leq p_N<\bar{p}^*.$$ 
    
    \item[-] In the sequel we will make use of the following numbers:
    \begin{equation}
    \label{parameters}
    \lambda = N(\bar{p}-2)+\bar{p}, \quad \alpha = \frac{N}{\lambda}, \quad \alpha_i = \frac{N(\bar{p}-p_i)+\bar{p}}{\lambda p_i}.    
    \end{equation}

    \item[-] For any $M, \rho>0$, the {\it intrinsic cube} and the {\it backward intrinsic cylinder} are defined respectively as 
    \[\K_{\rho}(M)= \prod_{i=1}^N \bigg{ \{} |x_i| < M^{\frac{p_i-2}{p_i}} \rho^{\frac{\bar{p}}{p_i}}\bigg{\}},\qquad \quad 
    \Q_{\rho}^{-}(M) = \prod_{i=1}^N \bigg{ \{} |x_i| < M^{\frac{p_i-2}{p_i}} \rho^{\frac{\bar{p}}{p_i}}\bigg{\}} \times \bigg( - \rho^{\bar{p}},0 \bigg].\]
    The notation of forward and full intrinsic cylinders is analogous to the one above.
    
    \item[-] The function $\pi_i:\R^N\to\R$, $\pi_i(x)=x_i$, $i=1,\ldots,N$, will denote the projection with respect to the $i$-th space variable. Moreover, $\pi:\R^N\times\R\to\R^N$, $\pi(x,t)=x$, stands for the projection in the space variables.
    
    \item[-] We denote by $\gamma$ a positive constant (depending only on the data, i.e., $N$ and $p_i$s) that may vary from line to line.

\end{itemize}
}

\section{Preliminaries} \label{Preliminaries}
We introduce the parabolic anisotropic spaces, which are the natural setting to work within. We define
\[ W^{1,{\bf{p}}}_o(\Omega):= \{ u \in W^{1,1}_o(\Omega) |\,  \partial_i u \in L^{p_i}(\Omega) \}, \]
\[ W^{1,{\bf{p}}}_{loc}(\Omega):= \{ u \in W^{1,1}_{loc}(\Omega) |\,  \partial_i u \in L^{p_i}(\Omega) \}, \]
% \[W^{1,{\bf{p}}}_{loc}(\Omega):= \{ u \in L^1_{loc}(\Omega) |\,  D_i u \in L^{p_i}_{loc}(\Omega) \} \]
% \[ L^{{\bf{p}}}(0,T;W^{1,{\bf{p}}}_o(\Omega)):= \{u \in L^1(0,T;W^{1,1}_o(\Omega))|\, D_iu \in L^{p_i}(0,T;L^{p_i}_{loc}(\Omega))   \} \]
\[ L^{{\bf{p}}}_{loc}(0,T;W^{1,{\bf{p}}}_o(\Omega)):= \{u \in L^1_{loc}(0,T;W^{1,1}_o(\Omega))|\, \partial_i u \in L^{p_i}_{loc}(0,T;L^{p_i}_{loc}(\Omega))   \}. \]
%We use the symbol $u_{x_i}$ meaning the $i$-th weak derivative for a function $u \in W^{1,\bf{p}}_0(\Omega)$, and the partial derivatives $\frac{\partial u}{\partial x_i}$ for formal calculations.\newline
A function \[ u \in C^0_{loc}(0,T; L^2_{loc}(\Omega)) \cap L^{\bf{p}}_{loc}(0,T;W^{1,{\bf{p}}}_{loc}(\Omega))\] is called a {\it local weak solution} of \eqref{EQ} if, for any $0<t_1<t_2<T$ and any compact set $K \subset \subset \Omega$, it satisfies
\begin{equation} \label{anisotropic-localweaksolution}
\int_{K} u \varphi \, dx \bigg|_{t_1}^{t_2}+ \int_{t_1}^{t_2} \int_{K} (-u \, \partial_t \varphi + \sum_{i=1}^N \, |\partial_i u|^{p_i-2} \partial_i u \, \partial_i \varphi) \, dx dt=0,
\end{equation} \noindent for all $\varphi \in C^{\infty}_{loc}(0,T;C_o^{\infty}(\Omega))$. By a density and approximation argument, we can consider test functions in \eqref{anisotropic-localweaksolution} in the bigger space \[ \varphi \in W^{1,2}_{loc}(0,T;L^2_{loc}(\Omega))\cap L^{\bf{p}}_{loc}(0,T;W^{1,{\bf{p}}}_o(\Omega)),\] provided $\Omega \subset \subset \R^N$ is a rectangular domain (see \cite{Haskovec-Schmeiser} for an extension to more general domains).

\subsection{Scaling properties of solutions} \label{scaling}
In the present subsection we show some important scaling properties of solutions to \eqref{EQ} and their correspondence with stationary solutions to a Fokker-Planck-type equation.
\vskip0.2cm 
\noindent 

\begin{proposition} \label{transformation-group}
Let $u$ be a local weak solution to the equation \eqref{EQ} in $\Omega_T$. For any $M,\rho>0$ appropriate for the inclusion $\Q_{\rho}(M) \subset \Omega_T$, we define the parametric transformation
\begin{equation} \label{Trans}
\tilde{T}_{\rho,M} (x,t)= \bigg(M^{\frac{p_i-2}{p_i}} \rho^{\frac{\bar{p}}{p_i}} x_i ,  \rho^{\bar{p}} {t}\bigg).
\end{equation} Then the transformed function
\begin{equation} \label{transf}
    \T(u)(x,t)= M^{-1} u \bigg( \tilde{T}_{\rho,M} (x,t) \bigg)= M^{-1}u \bigg(M^{\frac{p_i-2}{p_i}} \rho^{\frac{\bar{p}}{p_i}} x_i , \rho^{\bar{p}} {t}\bigg)
\end{equation} \noindent is a solution to \eqref{EQ} in $\tilde{T}_{\rho,M}^{-1}(\Omega_T)$. 
\end{proposition}

\begin{proof} We perform some formal algebraic computations, representing change of variables in the integrals of definition \eqref{anisotropic-localweaksolution}.\vskip0.2cm 

\noindent If we generally suppose
\[
\T (u) = M^{-1}u\bigg( L_i x_i,Tt \bigg)
\]
for some $M,L_i,T>0$, then
\[
\partial_t \bigg( \T u \bigg)= M^{-1}T \partial_t u (L_i x_i,Tt)\quad \mbox{and} \quad \partial_i \bigg( \T u \bigg)= M^{-1}L_i  \partial_i u(L_i x_i,Tt).
\] Thus, imposing the equation for $\T u$, namely,
\[
\partial_t \bigg( \T u \bigg) = \sum_{i=1}^N \partial_i \bigg( | \partial_i (\T u) |^{p_i-2} \partial_i (\T u) \bigg),
\]
we find
\[
\begin{split}
M^{-1}T \partial_t u(L_i x_i,Tt) = \partial_t \bigg( (\T u)(x,t) \bigg) &= \sum_{i=1}^N \partial_i \bigg( | \partial_i (\T u)(x,t) |^{p_i-2} \partial_i (\T u)(x,t) \bigg)\\
&= \sum_{i=1}^N L_i^{p_i}M^{1-p_i} \partial_i \bigg( |\partial_i u |^{p_i-2} \partial_i u \bigg)(L_i x_i,Tt).
\end{split}
\]
%We write for some $A>0$ to be defined
%\[
%L_i= M^{\frac{p_i-\bar{p}}{p_i}} A^{\frac{1}{p_i}},
%\] to restore the homogeneity in the equation. In order to let $\T %u $ satisfy the same equation we need
%\[
%M^{-1}T=M^{1-\bar{p}}A, \quad \Rightarrow \quad A=TM^{\bar{p}-2}, 
%\] and so the transformation is
Furthermore, we impose
\[
M^{-1}T = L_i^{p_i}M^{1-p_i} \quad \forall i=1,\ldots,N
\]
to restore the homogeneity in the equation. We find $L_i=[M^{p_i-2}T]^{\frac{1}{p_i}}$, whence
\begin{equation}
\label{gentransf}
\T u = M^{-1} u \bigg( \bigg[ M^{p_i-2}T \bigg]^{\frac{1}{p_i}} x_i,Tt\bigg).
\end{equation}
Taking $T=\rho^{\bar{p}}$ concludes the proof.\end{proof} \noindent 

\begin{remark} The peculiarity of the scaling \eqref{transf} is that it does not alter the time variable, from the point of view of intrinsic geometry. Indeed, the parameter $M$ is usually chosen to be a suitable multiple of either the oscillation or the $L^{\infty}$ norm of the solution itself, therefore leading to a geometry within the equation evolves in an intrinsic fashion (see \cite{DB},\cite{Urb}). \newline 
Moreover, the proof of Proposition \ref{transformation-group} reveals that \eqref{transf} is not the only invariant: we may consider, for instance, also the transformation
\begin{equation} \label{scalings}
\T_{\rho,M} u= M^{-1} u\bigg( M^{\frac{p_i-\bar{p}}{p_i}}\rho^{\frac{\bar{p}}{p_i}} x_i, \, M^{2-\bar{p}}\rho^{\bar{p}}t\bigg),
\end{equation}
corresponding to $T=M^{2-\bar{p}} \rho^{\bar{p}}$ in \eqref{gentransf}.
This transformation has been used extensively in \cite{Ciani-Mosconi-Vespri}, with the aim of obtaining a Harnack inequality which intrinsically scales within the particular geometry dictated by the transformation.
\end{remark}

\begin{definition}\label{basics} We define the intrinsic anisotropic cube  by  transformation \eqref{Trans} on the space variables, \begin{equation}\label{intrinsic-anisotropic-cube}
    \K_{\rho}(M)= \prod_{i=1}^N \bigg{ \{} |x_i| < M^{\frac{p_i-2}{p_i}} \rho^{\frac{\bar{p}}{p_i}}\bigg{\}},
\end{equation} \noindent and the intrinsic anisotropic cylinders 
\begin{equation} \label{intrinsic-cylinders}
    \Q_{\rho}^{-}(M):= \tilde{T}_{\rho,M} (Q_1^{-}) = \prod_{i=1}^N \bigg{ \{} |x_i| < M^{\frac{p_i-2}{p_i}} \rho^{\frac{\bar{p}}{p_i}}\bigg{\}} \times \bigg( - \rho^{\bar{p}},0 \bigg].
\end{equation} \noindent Similarly we define forward and full intrinsic cylinders.
\end{definition} 

\begin{remark}
\label{standardcase}
Definition \ref{basics} is motivated by Proposition \ref{transformation-group} and leads to the following consequence. \newline
If $u$ solves \eqref{EQ} in $\Q_{\rho}^-(M)$, then $\T u$ solves \eqref{EQ} in $Q_1^{-}$. Vice-versa if $u$ solves \eqref{EQ} in $Q_{1}^-$ then $\T^{-1} (u)= Mu \bigg(M^{\frac{2-p_i}{p_i}}\rho^{-\frac{\bar{p}}{p_i}}x_i, \rho^{-\bar{p} }t \bigg)$ solves \eqref{EQ} in $\Q_{\rho}^-(M)$.
\end{remark}

\begin{proposition}
\label{preservenorm}
Let $u$ be a local weak solution to equation \eqref{EQ}. Then the parametric transformations preserving the $L^1$ norm of $u$ correspond to \eqref{transf} for $M=\rho^{-\alpha\bar{p}}$, that is,
\begin{equation} \label{transformation}
    \T_{\rho} u = \rho^{\alpha\bar{p}} u\bigg( \rho^{\alpha_i\bar{p}} x_i,\rho^{\bar{p}} t \bigg),
\end{equation}
where $\alpha,\alpha_i$ were defined in \eqref{parameters}.
\end{proposition}
\begin{proof} Performing a change of variables, besides recalling \eqref{parameters}, we get
\begin{equation} \label{c-var}
\int_{\pi(\tilde{T}_{\rho,M}(\Q_1))} \T u(x,t) dx= M^{-1}\prod_{i=1}^N \bigg(M^{\frac{p_i-2}{p_i}} \rho^{\frac{\bar{p}}{p_i}} \bigg)^{-1} \int_{K_1} u(y,s) dy =\bigg(M^{-\frac{\lambda}{\bar{p}}} \rho^{-N}  \bigg)\int_{K_1} u(y,s) dy. \end{equation}
Hence, imposing $M^{-\frac{\lambda}{\bar{p}}} \rho^{-N}=1$, we find $M=\rho^{-\alpha\bar{p}}$, as desired.
\end{proof} 
\begin{remark}
It is worth noticing the following important geometric property, used also in the proof of Proposition \ref{preservenorm}: for any $M,\rho>0$, the total volumes of the anisotropic cube and the anisotropic cylinder depend on $p_i$s, i.e., \[
|\K_{\rho}(M)|= 2^N \rho^N M^{\frac{N(\bar{p}-2)}{\bar{p}}}= \rho^{N} M^{\frac{N(\bar{p}-2)}{\bar{p}}}\, |K_1|, \qquad |\Q_{\rho}^-(M)|= 2^N \rho^{N+\bar{p}} M^{\frac{N(\bar{p}-2)}{\bar{p}}} = \rho^{N+\bar{p}} M^{\frac{N(\bar{p}-2)}{\bar{p}}}\, |Q_1^-|.
\] 
\end{remark}

\begin{definition}
A solution $u$ to \eqref{EQ} in $\R^{N+1}$ is said to be a {\it self-similar solution} if it satisfies $\T_{\rho} u =u$ for all $\rho>0$, where $T_\rho$ was defined in \eqref{transformation}.
\end{definition}

\vskip0.2cm \noindent Now we consider the continuous transformation $\Phi$ and its inverse $\Phi^{-1}$ defined as
\begin{equation}\label{continuous-transformation}
    \Phi (u)(x,t)= w(x,t)=e^{\alpha t} u( e^{\alpha_i t}x_i,e^t), \quad \Phi^{-1}(w)(y,s)=u(y,s)=s^{-\alpha} w(s^{-\alpha_i}y_i, \log s).
\end{equation}

\noindent 
This map formally sends solutions to \eqref{EQ} in $\Sigma^+:=  \R^N \times (0, +\infty) $ into solutions of the anisotropic Fokker-Planck-type equation
\begin{equation}\label{FokkerPlanck}
  \partial_t w = \sum_{i=1}^{N}\partial_i [ (| \partial_i w  |^{p_i-2} \partial_i w ) +\alpha_i y_i w ] \quad  \mbox{in} \quad \Sigma := \R^N \times \R.
%   \quad \alpha_i= \bigg(\frac{1+2\alpha}{p_i}-\alpha  \bigg),
\end{equation}
% where $\alpha= \frac{N}{N(\bar{p}-2)+\bar{p}}.$
\noindent For each fixed time $t=\log (\rho^{-\lambda})$, $\rho>0$, $\Phi$ corresponds to a parametric transformation of type \eqref{transformation}, thus preserving the $L^1$ norm; indeed, it is readily seen that
\begin{equation} \label{Phi}
\Phi(u)(x,\log(\rho^{-\lambda}))=\T_{\rho^{-\lambda/\bar{p}}}u(x,1). \end{equation} 

\noindent Now we present a characterization of the self-similar solutions.
\begin{proposition} \label{fk-ss}
Self-similar solutions to \eqref{EQ} in $\Sigma^+$ correspond to stationary solutions to the Fokker-Planck equation \eqref{FokkerPlanck} and vice-versa.
\end{proposition}
\begin{proof}
% A self similar solution of \eqref{EQ} which maintains its $L^1$-norm unvaried is a solution that satisfies $\T_{\rho} u=u$ for each $\rho>0$. 
\noindent Let us consider a self-similar solution $u$ to \eqref{EQ} in $\Sigma^+$. We already know that $w=\Phi u$ is a solution to \eqref{FokkerPlanck}. It remains to show that $w$ is stationary. By \eqref{Phi} and the self-similarity of $u$, for all $(x,t)\in\Sigma$ we get
\[
w(x,t) = w(x,\log(\rho^{-\lambda})) = \T_{\rho^{-\lambda/\bar{p}}}u(x,1) = u(x,1) = w(x,0),
\]
being $t=\log(\rho^{-\lambda})$ for some $\rho>0$.\newline
Vice-versa, let $w$ be a stationary solution to \eqref{FokkerPlanck}. We already know that $u$ solves \eqref{EQ} in $\Sigma^+$, so it suffices to show that $u$ is self-similar. For any $\rho>0$, we choose $t=\log(\rho^{\bar{p}}l)$, $l>0$, in \eqref{continuous-transformation} and use the fact that $w$ is stationary to obtain
\[
l^\alpha \rho^{\alpha\bar{p}}u(\rho^{\alpha_i\bar{p}}l^{\alpha_i} x_i,\rho^{\bar{p}}l) = w(x,t) = w(x,\log l) = l^\alpha u(l^{\alpha_i}x_i,l).
\]
Dividing by $l^\alpha$, besides performing the change of variables $y_i=l^{\alpha_i}x_i$, leads to
\[
\T_\rho u = u \quad \forall \rho>0,
\]
which is the self-similarity of $u$.
\end{proof} \noindent

\begin{definition}\label{Barenblatt} A self-similar solution to \eqref{EQ} in $\Sigma^+$ (or, equivalently, a solution to \eqref{EQ} corresponding to a stationary solution to the Fokker-Planck equation \eqref{FokkerPlanck}) is said to be a {\it{Barenblatt Fundamental solution}}; it is denoted by $\B$, in analogy with the literature regarding the $p$-Laplacian\footnote{Indeed, the epithet {\it Fundamental} does not mean that solutions are represented by an integral convolution with kernel $\B$, but that the classic $\B$ function approaches to the heat kernel as $p \rightarrow 2$. The Barenblatt solution for the p-Laplacian equation can be found in \cite{Barenblatt}.}.
\end{definition}

\section{Scaling invariants}
\label{invariants} 
Definition \ref{Barenblatt} is invariant under the scalings \eqref{transf} and \eqref{scalings}. In this section we show that also the energy of solutions is invariant, and the same holds for a particular energetic property of solutions, that can be regarded as a measure-theoretical maximum principle.

\begin{lemma}[Energy Estimates] \label{EE}
Let $u$ be a local weak solution to equation \eqref{EQ} in $\Omega_T$. Let $(x_o,t_o) \in \Omega_T$ and $\rho,M>0$ be such that $(x_o,t_o)+\Q_{\rho}^-(M)\subset \Omega_T$. Then, for each function of the form \[ C_o^{\infty}((x_o,t_o)+\Q_{\rho}^-(M)) \ni \eta = \prod_{i=1}^N \eta_i^{p_i}(x_i,t) \quad \text{with} \quad \eta_i \in C_o^{\infty}(\pi_i(x_o+\K_\rho(M)) \times (t_o-\rho^{\bar{p}},t_o]),\] we have the following estimates, valid for all $t_o-\rho^{\bar{p}}<s<t<t_o$ and $k \in \R$:
\begin{equation} \label{energy}
\begin{aligned}
    \int_{\K_{\rho}(M)}& (u-k)_{\pm}^2 \eta(x,\tau) dx \bigg|_{\tau=s}^{\tau=t} + \sum_{i=1}^N \int_{s}^t \int_{\K_{\rho}(M)} |\partial_i [\eta (u-k)_{\pm}]|^{p_i}  \, dx d\tau \\
    & \leq  \gamma \bigg{\{} \int_{s}^t \int_{\K_{\rho}(M)} (u-k)_{\pm}^2 \partial_{\tau}\eta(x,\tau) \, dx d\tau + \sum_{i=1}^N \int_{s}^t \int_{\K_{\rho}(M)} |(u-k)_{\pm}|^{p_i} \hat{\eta}_i |\partial_i \eta_i|^{p_i} \, dx d\tau  \bigg{\}},
    \end{aligned}
\end{equation}
where $\hat{\eta}_i:= \eta/\eta_i^{p_i}$ and $\gamma>0$ is a suitable constant (depending only on $N$ and $p_i$s).
\noindent

\end{lemma}
\begin{proof}
The function $u$ solves equation \eqref{EQ} in $(x_o,t_o)+ \Q_{\rho}^-(M)$, so $\T (u)$ (defined in \eqref{transf}) solves \eqref{EQ} in $Q_1^-$, according to Remark \ref{standardcase}. Now, Lemma 3.1 of \cite{Mosconi} on unitary cylinders ensures that for each function of the form \begin{equation} \label{test}
C_o^{\infty}(Q_1) \ni \eta = \prod_{i=1}^N \eta_i^{p_i}(x_i,t) \quad \text{with} \quad \eta_i \in C_o^{\infty}(\pi_i(K_1) \times (-1,0]),
\end{equation} \noindent 
we have, for all $-1<s_1<s_2<0$ and $\bar{k}\in\R$,
\begin{equation} \label{EEE}
\begin{aligned}
\int_{K_1} (\T u -\bar{k})_{\pm}^2& \eta \, dy \bigg|_{s_1}^{s_2}+ \sum_{i=1}^N \int_{s_1}^{s_2} \int_{K_1} | \partial_i [\eta(\T u-\bar{k})_{\pm}]|^{p_i} \, dy ds \\
&\leq \gamma \int_{s_1}^{s_2} \int_{K_1} |(\T u-\bar{k})_{\pm} |^2 \partial_s \eta \, dyds + \sum_{i=1}^N \int_{s_1}^{s_2} \int_{K_1} |(\T u- \bar{k})_{\pm} |^{p_i} \hat{\eta}_i | \partial_i \eta_i |^{p_i} \,dyds.
\end{aligned}
\end{equation}

%\begin{equation} \label{ei}
%\begin{aligned}
%   \int_{K_1}& (\T u-k)_{\pm}^2 \eta(x,\tau) dx \bigg|_{s}^t + \sum_{i=1}^N \int_{s}^t \int_{K_1} |\partial_i[\eta  (\T u-k)_{\pm}]|^{p_i} \, dx d\tau \\
%  & \leq  \gamma \bigg{\{} \int_{s}^t \int_{K_1} (\T u-k)_{\pm}^2 \partial_{\tau}\eta(x,\tau) \, dx d\tau + \sum_{i=1}^N \int_{s}^t \int_{K_1} |(\T u-k)_{\pm}|^{p_i} \hat{\eta}_i |\partial_i \eta_i|^{p_i} \, dx d\tau  \bigg{\}}.
% \end{aligned}
%\end{equation} \noindent

\noindent Now we show that \eqref{energy} comes from \eqref{EEE} by performing the change of variables \eqref{Trans}-\eqref{transf}, besides observing that $[\T u >\bar{k}]\cap K_1=[u>k]\cap \K_{\rho}(M)$ provided $\bar{k}=k/M$. Indeed, let us consider the change of variables 
\[u(x,t) = u(M^{(p_i-2)/p_i} \rho^{\bar{p}/p_i}y_i,\rho^{\bar{p}}s) = M\T(u)(y,s),
\]
with the stipulations
\[
\begin{cases}
x_i= M^{\frac{p_i-2}{p_i}}\rho^{\frac{\bar{p}}{p_i}}y_i,\\
t= \rho^{\bar{p}} s.
\end{cases}
\]
We observe that
\[
\begin{cases}
\partial_{x_i} u (x,t)= M^{2/p_i}\rho^{-\bar{p}/p_i}\partial_{y_i} \T u (y,s),\\
\partial_{t} u( x,t)= M \rho^{-\bar{p}} \partial_s\T u (y,s),
\end{cases}
\]
and $dx(y)=\left(\prod_i |dx_i /dy_i|\right) dy =\rho^N M^{\frac{N(\bar{p}-2)}{\bar{p}}} dy$. Hence the first integral in \eqref{EEE} becomes
\begin{equation*}
    \begin{aligned}
 \int_{K_1} (\T u -\bar{k})_{\pm}^2 \eta \, dy \bigg|_{s_1}^{s_2}=& \int_{K_{\rho}(M)} (M^{-1}(u(x,t)-k))_{\pm}^2 \eta(y(x),s(t)) \, (\rho^{-N} M^{\frac{-N(\bar{p}-2)}{\bar{p}}} dx)\bigg|_{t_1}^{t_2}\\
 &=  \rho^{-N}M^{-\frac{[N(\bar{p}-2)+2\bar{p}]}{\bar{p}}} \int_{K_{\rho}(M)} (u(y,t)-k)_{\pm}^2\eta \, dx \bigg|_{t_1}^{t_2},
    \end{aligned}
\end{equation*} being $t_1:=\rho^{\bar{p}}s_1<\rho^{\bar{p}}s_2=:t_2$. Similarly we evaluate the other integrals of \eqref{EEE}, obtaining
\begin{equation*}
    \begin{aligned}
\int_{s_1}^{s_2} \int_{K_1} | \partial_{y_i} [\eta(\T u-\bar{k})_{\pm}]|^{p_i} \, dy ds=&
% \int_{t_1}^{t_2} \int_{K_{\rho}(M)} |(M/\rho)^{-\bar{p}/p_i} \partial_{x_i}( u-\bar{k}/M)_{\pm}|^{p_i} \eta (\rho^{-N} dx)(M^{\bar{p}-2}\rho^{-\bar{p}} dt)\\ &=
 \rho^{-N}M^{-\frac{[N(\bar{p}-2)+2\bar{p}]}{\bar{p}}} \int_{t_1}^{t_2} \int_{K_{\rho}(M)} |\partial_{x_i} [\eta(u-k)_{\pm}]|^{p_i} \, dxdt.\end{aligned}
\end{equation*}

\begin{equation*}
    \begin{aligned}
\int_{s_1}^{s_2} \int_{K_1} |(\T u-\bar{k})_{\pm} |^2 \partial_s \eta \, dyds=&
% \int_{t_1}^{t_2} \int_{K_{\rho}(M)} |(M^{-1} u(y,t)-\bar{k})_{\pm}^2  (M^{2-\bar{p}} \rho^{\bar{p}}) \partial_t \eta \, (\rho^{-N})dx(M^{\bar{p}-2}\rho^{-\bar{p}})dt\\&=
\rho^{-N}M^{-\frac{[N(\bar{p}-2)+2\bar{p}]}{\bar{p}}} \int_{t_1}^{t_2} \int_{K_{\rho}(M)} |(u(y,t)-k)_{\pm}^2 \partial_t \eta \, dxdt.
    \end{aligned}
\end{equation*}
\begin{equation*}
    \begin{aligned}
\int_{s_1}^{s_2} \int_{K_1} |(\T u- \bar{k})_{\pm} |^{p_i} \hat{\eta}_i | \partial_{y_i} \eta_i |^{p_i} \,dyds=&
% \int_{t_1}^{t_2} \int_{K_{\rho}(M)} 
% |(M^{-1} u(x,t)-\bar{k})_{\pm}|^{p_i} \hat{\eta} | (M^{\frac{\bar{p}-p_i}{p_i}}\rho^{-\frac{\bar{p}}{p_i}}) \partial_{x_i} \eta_i|^{p_i} \, \rho^{-N} dx M^{\bar{p}-2} \rho^{-\bar{p}} dt =\\&
\rho^{-N}M^{-\frac{[N(\bar{p}-2)+2\bar{p}]}{\bar{p}}} \int_{t_1}^{t_2} \int_{K_{\rho}(M)} 
|u(x,t)-k)_{\pm}|^{p_i} \hat{\eta}_i | \partial_{x_i} \eta_i|^{p_i} \, dxdt.
    \end{aligned}
\end{equation*}
Collecting the terms $\rho^{-N}M^{-\frac{[N(\bar{p}-2)+2\bar{p}]}{\bar{p}}}$ we get \eqref{energy}. Hence, energy estimates are invariant under the scaling transformation \eqref{transf}.
\end{proof}
\begin{remark}
Clearly, the energy estimates above are valid also in forward and full cylinders $(x_o,t_o)+\Q_{\rho}^+(M)$, $(x_o,t_o)+\Q_{\rho}(M)$, provided they are contained in $\Omega_T$.\end{remark}

\noindent The next Lemma is a sort of measure-theoretical maximum principle, popular amongst nonlinear analysts as Critical Mass Lemma (following Caffarelli), or De Giorgi-type Lemma (following DiBenedetto). It may be proven at ease for unitary cylinders, and then re-interpreted in the intrinsic geometry dictated by the scaling \eqref{transf}. To show the convenience of using \eqref{transf}, first we prove the lemma in its general form, and then we discuss its invariance with respect to the scaling. \newline
We recall that local weak sub-solutions (resp, super-solutions) to \eqref{EQ} are locally bounded from above (resp., below) in $\Omega_T$ (see, e.g., \cite{Mosconi}, \cite{Mingqi}), provided an additional condition constraining the spareness of $p_i$s is ensured. Let us fix a cylinder $(y,s)+\Q_{2\rho}(\theta) \subset \subset  \Omega_T$, being $(y,s)\in\Omega_T$ and $\rho,\theta>0$ appropriate. Let $\mu^+, \mu^-$ be such that
\[\mu^- \leq \essinf_{(y,s)+\Q_{2\rho}(\theta)} u\leq \essup_{(y,s)+\Q_{2\rho}(\theta)} u \leq \mu^+.\] We also fix $\omega>0$, $\xi \in (0,1]$, and $a \in (0,1)$.

\begin{lemma}[De Giorgi-type/Critical Mass]\label{DG}
Let $u$ be a local weak super-solution to \eqref{EQ} in $\Omega_T$ locally bounded from below, and let $\rho,\theta,\mu^{\pm},\omega,\xi,a$ be defined as above. Then there exists $\nu^- \in (0,1)$, depending on the data $N$,$p_i$s and on the parameters $\theta,\omega,\xi,a$ but not on the radius $\rho$, such that if \begin{equation} \label{DG-HP}
|[u\leq \mu^- + \xi \omega] \cap [(y,s)+ \Q_{2\rho}^-(\theta)] | \leq \nu^- |\Q_{2\rho}^-(\theta)|
\end{equation} \noindent then 
\begin{equation} \label{DG-TH}
    u \ge \mu^-+a \xi\omega \quad \text{a.e. in} \, \, \Q_{\rho}^-(\theta).
\end{equation}\noindent Likewise, if $u$ is a local weak sub-solution to \eqref{EQ} in $\Omega_T$ which is locally bounded from above, then there exists $\nu^+\in (0,1)$, depending on the data $N$,$p_i$s and on the parameters $\theta,\omega,\xi,a$ but not on the radius $\rho$, such that if \begin{equation} \label{DG+HP}
|[u\ge  \mu^+ - \xi \omega] \cap [(y,s)+ \Q_{2\rho}^-(\theta)] | \leq \nu^+ |\Q_{2\rho}^-(\theta)|
\end{equation} \noindent then 
\begin{equation} \label{DG+TH}
    u \leq \mu^+-a \xi\omega \quad \text{a.e. in} \, \, \Q_{\rho}^-(\theta).
\end{equation}\noindent 

\end{lemma}

\begin{proof}
We prove \eqref{DG-TH}, since the proof of \eqref{DG+TH} is analogous. Without loss of generality we assume $(y,s)=(0,0)$, just to ease the notation. Let us set, for any $n, \bar{n}\in\N\cup\{0\}$,
\[\rho_n=\rho+\frac{\rho}{2^n}, \qquad \K_{n}=  \prod_{i=1}^N\bigg{\{}|x_i|< \theta^{\frac{p_i-2}{p_i}} \rho^{\frac{\bar{p}}{p_i}}\bigg(1+\frac{1}{2^{n+\bar{n}}} \bigg)\bigg{\}}, \qquad \Q_n=\K_n \times (-\rho^{\bar{p}}_n, 0].\]
Since $K_0 \to \K_\rho(\theta)$ as $\bar{n}\to\infty$, we fix $\bar{n}$ such that $\K_0 \subset \K_{2\rho}(\theta)$. Notice that $\bar{n}$ can be chosen in such a way that it depends only on $N$ and $p_i$s.
%, as for instance \[ 
%\bar{n}= [\, |\log_2(2^{\bar{p}/p_N}-1))|\, ]^+\] will suffice, where $[\cdot]^+$ denotes the upper integer part. Then we consider $\bar{n}$ fixed and we proceed by admitting that constant $\gamma$ may depend also on $\bar{n}$. 
We apply energy estimates \eqref{energy} over $\K_n, \Q_n$ to the truncations $(u-k_n)_-$ at the levels
\[k_n= \mu^-+ \xi_n \omega, \quad \text{where} \quad \xi_n= a \xi + \frac{(1-a) \xi}{2^n}. \]
Incidentally, notice that
\[
|(u-k_n)_-|\leq \xi_n \omega\leq \xi \omega.
\]

\noindent For any $n$, we pick a cut-off function $\eta_n$ of the form $\eta_n= \bar{\eta}(t) \prod_{i=1}^N \eta_i^{p_i}(x_i)$, where

\begin{equation*}
    \eta_i(x_i) =\begin{cases}
    1, \quad \text{in} \quad \pi_i(\K_{n+1}),\\
    0, \quad \text{in} \quad \R \setminus \pi_i(\K_n),
    \end{cases} \quad |\partial_i \eta_i|\leq \frac{\gamma 2^{n}}{\theta^{\frac{p_i-2}{p_i}}\rho^{\frac{\bar{p}}{p_i}}}. 
\end{equation*}

\begin{equation*}
    \bar{\eta}(t)=\begin{cases}
    1, \quad \text{when} \quad t\ge -\rho_{n+1}^{\bar{p}},\\
    0, \quad \text{when} \quad t<-\rho^{\bar{p}}_{n},
    \end{cases} \quad  |\partial_t \bar{\eta}|\leq \frac{\gamma}{(\rho_n^{\bar{p}}- \rho_{n+1}^{\bar{p}})} \leq \frac{\gamma 2^{(n+1)\bar{p}}}{\rho^{\bar{p}}},
\end{equation*}
The energy estimate \eqref{energy}, applied with these choices, yields

\begin{equation*}
\begin{aligned}
\E_n:= \sup_{(-\rho_n^{\bar{p}},0]} \int_{\K_n} (u-k_n)_-^2 \eta_n \, dx +& \sum_{i=1}^N \int \int_{\Q_n} |\partial_i [\eta_n(u-k_n)_-]|^{p_i} \, dxdt  \\
&\leq\gamma \int \int_{\Q_n}\bigg{\{} (u-k_n)_-^2 |\partial_t \bar{\eta}|+\sum_{i=1}^N |(u-k_n)_-|^{p_i} |\partial_i \eta_{i}|^{p_i}  \bigg{\}}dxdt  \\
& \leq\frac{\gamma 2^{\bar{p}(n+1)}}{\rho^{\bar{p}}} (\xi \omega)^2 \bigg(1+\sum_{i=1}^N \bigg(\frac{\xi \omega}{\theta} \bigg)^{p_i-2}  \bigg)|[u< k_n]\cap \Q_n|.
\end{aligned}
\end{equation*}
Now we combine these estimates of the energies $\E_n$ with the embedding inequality (see \cite{Mosconi})

\begin{equation*}\label{embedding}
\begin{aligned}
    \int \int_{\Q_n}& |(u-k_n)_- \eta_n |^{\bar{p}(\frac{N+2}{N})}dxdt\\
    & \leq    \gamma \bigg( \sup_{(-\rho^{\bar{p}}_n,0]} \int_{\K_n} (u-k_n)_-^2 \eta_n^2\, dx+ \sum_{i=1}^N \int \int_{\Q_n} |\partial_i[\eta_n [(u-k_n)_-] |^{p_i}\, dxdt  \bigg)^{\frac{N+\bar{p}}{N}}= \gamma \E_n^{\frac{N+\bar{p}}{N}}.
    \end{aligned}
\end{equation*} \noindent Observing that  $(u-k_n)_-\ge (k_n-k_{n+1})= (1-a)\xi \omega/ 2^{n+1}$ in $[u<k_{n+1}]\cap \Q_{n+1}$, we get the chain

\begin{equation}\label{chain}
\begin{aligned}
\bigg(\frac{(1-a)\xi \omega}{2^{n+1}}& \bigg)^{\bar{p} (\frac{N+2}{N})} |[u<k_{n+1}] \cap \Q_{n+1}| \leq \int \int_{\Q_{n+1}} |(u-k_n)_- |^{\bar{p}(\frac{N+2}{N})}dxdt \\
&\leq \int \int_{\Q_n} |(u-k_n)_- \eta_n |^{\bar{p}(\frac{N+2}{N})}dxdt \leq \gamma \E_n^{\frac{N+\bar{p}}{N}}\\
& \leq\gamma \bigg[ \frac{2^{\bar{p}(n+1)}}{\rho^{\bar{p}}}(\xi \omega)^2 \bigg(1+\sum_{i=1}^N \bigg( \frac{\xi \omega}{\theta} \bigg)^{p_i-2}  \bigg) \bigg]^{\frac{N+\bar{p}}{N}} |[u<k_n]\cap \Q_n|^{\frac{N+\bar{p}}{N}}.
\end{aligned}
\end{equation}
We set $Y_n=\frac{|[u<k_n]\cap \Q_n|}{|\Q_n|}$ and recall that up to a constant we have  $|\Q_n|^{\frac{\bar{p}}{N}}\leq \theta^{\bar{p}-2}\rho^{\frac{\bar{p}}{N}(\bar{p}+N)}$. Dividing \eqref{chain} for $|\Q_{n+1}|$, as well as noticing that $|\Q_{n+1}|\ge 2^{-(N+1)} |\Q_n|$, we obtain
\begin{equation}
\label{iteration}
\begin{aligned}
Y_{n+1}\leq& \gamma \bigg(\frac{2^{n+1}}{(1-a)\xi \omega} \bigg)^{\bar{p} (\frac{N+2}{N})} \bigg[ \frac{2^{\bar{p}(n+1)}}{\rho^{\bar{p}}}(\xi \omega)^2 \bigg(1+\sum_{i=1}^N \bigg( \frac{\xi \omega}{\theta} \bigg)^{p_i-2}  \bigg) \bigg]^{\frac{N+\bar{p}}{N}} Y_{n}^{1+\frac{\bar{p}}{N}} |\Q_n|^{\frac{\bar{p}}{N}}\\
&\leq 2^{n(2N+\bar{p}+2)\frac{\bar{p}}{N}} \bigg[ \frac{\gamma (\xi \omega/ \theta)^{2-\bar{p}}}{(1-a)^{\frac{\bar{p}}{N}(N+2)}} \bigg(1+ \sum_{i=1}^N (\xi \omega/\theta)^{p_i-2}  \bigg)^{\frac{N+\bar{p}}{N}}  \bigg] Y_n^{1+\frac{\bar{p}}{N}}. \end{aligned}
\end{equation}

\noindent 
According to \cite[Lemma 4.1, page 12]{DB}, \eqref{iteration} produces $Y_\infty:=\lim_{n\rightarrow \infty }Y_{n}=0$ provided 
\begin{equation}
\label{iterstart}
Y_{0}\leq 2^{-\frac{N}{\bar{p}}(2N+\bar{p}+2)}\bigg[ \frac{\gamma (\xi \omega/ \theta)^{2-\bar{p}}}{(1-a)^{\frac{\bar{p}}{N}(N+2)}} \bigg(1+ \sum_{i=1}^N (\xi \omega/\theta)^{p_i-2}  \bigg)^{\frac{N+\bar{p}}{N}}  \bigg]^{-\frac{N}{\bar{p}}}:= \nu^*.
\end{equation}
Then $Y_{\infty}=0$ implies \eqref{DG-TH}, concluding the proof. It remains to ensure \eqref{iterstart}. \newline
We set $\nu^-:=\nu^*/\gamma$, where $\nu^*$ stems from \eqref{iterstart} and $\gamma$ is such that $|\Q_{2\rho}^-(\theta)| \leq \gamma |\Q_{0}|$. Therefore 
\[|[u< \mu^-+\xi \omega] \cap \Q_{0}|\leq |[u< \mu^-+\xi \omega] \cap \Q_{2\rho}^-(\theta)|\leq \nu^- \gamma |\Q_0| = \nu^* |\Q_{0}|, \]
which guarantees \eqref{iterstart}. 

\noindent
In order to prove \eqref{DG+TH}, we may proceed as above, considering energy estimates \eqref{energy} in the same iterative geometry, but this time using the truncations $(u-k_n)_+$, being $k_n= \mu^+-\xi_n\omega$.
\end{proof} 

\begin{remark}
We observe that setting $\theta= \xi \omega$ frees $\nu^-$ from any dependence on $\theta,\omega,\xi$. When $\omega\ge\mu^+-\mu^-$, this choice of $\theta$ represents the intrinsic geometry, since $\omega$ is related to the oscillation of the solution $u$ in $\Q_{2\rho}^-(\theta)$.
\end{remark}

\noindent Simple adjustments on $\Q_n$ guarantee the validity of Lemma \ref{DG} also for forward cylinders and full cylinders.\vskip0.2cm 

\noindent Lemma \ref{DG} is invariant under the particular scalings \eqref{transf} and \eqref{scalings}. To show this, one can use the change of variables employed in the proof of Lemma \ref{energy} to prove the following lemma, and then one can go back to Lemma \ref{DG}. Accordingly, the argument sketched here shows the equivalence of Lemmas \ref{DG} and \ref{De Giorgi}.

\begin{lemma}\label{De Giorgi} Let $u$ be a local weak super-solution to \eqref{EQ} in $Q_1$ such that $\essinf_{Q_1^-} u \ge 0$. Then for any $a\in (0,1)$ there exists $\nu_a^+>0$, depending on $a,p_i,N$ but neither on $u$ nor on $\rho$, such that 
\begin{equation}\label{DGmeasurecondition}
        |[u\leq  1]\cap \Q_1^-(1)| \leq \nu_a^+ |\Q_1^-(1)| \quad \Rightarrow \quad \essinf_{\Q_{1/2}^-(1)} u \ge a \,. \end{equation}
%  being \[\Q_{1/2}^-=\prod_{i=1}^N \{|x_i| < (1/2)^{\frac{\bar{p}}{p_i}} \} \times (-(1/2)^{\bar{p}}, 0].\] 
Let $u$ be a local weak sub-solution to \eqref{EQ} in $Q_1$ such that $\essup_{Q_1^-} u \leq 1$. Then for any $a\in (0,1)$ there exists $\nu_a^->0$, depending on $a,p_i,N$ but neither on $u$ nor on $\rho$, such that
\begin{equation}\label{DG2}
        |[u \ge 1/2]\cap \Q_1^-(1/2)| \leq \nu_a^- |\Q_1^-(1/2)| \qquad \Rightarrow \qquad  \essup_{\Q_{1/2}^-(1/2)} u \leq (1-a/2) \,.\end{equation}
\noindent \end{lemma}
\begin{proof}
It suffices to apply Lemma \ref{DG} to the function $(u-\mu^-)/(\xi \omega)$ (resp., $(\mu^+-u)/(\xi \omega)$) with the choices $\mu^-=0$, $\theta=\xi\omega=1$ (resp., $\mu^+=1$, $\theta=\xi \omega=1/2$), and $\rho=1$ and  in the first (resp., second) case.
\end{proof}

\section{A topological consequence of energy invariants: lower semi-continuity of super-solutions}
\label{Ssemicontinuity}

\begin{theorem}\label{semicontinuity}
Let $u$ be a weak local super-solution to \eqref{EQ} in $\Omega_T$ locally bounded from below. Then $u$ is lower semi-continuous.
\end{theorem}

\begin{proof} We proceed in a way reminiscent of \cite{Liao}. Set $\Q_\rho := \Q_\rho(1) $ for all $\rho>0$, and consider the lower semi-continuous regularization of $u$, defined as 
\begin{equation}\label{lsc}
    u_*(x,t)= \lim_{\rho \rightarrow 0^+} \essinf_{(x,t)+\Q_\rho}\,  u \quad \forall (x,t) \in \Omega_T.
\end{equation} We observe that this function is well defined, since $(x,t)+\Q_\rho\subset\Omega_T$ for small values of $\rho$. It is a well-known fact that $u_*$ is lower semi-continuous. Accordingly, proving that $u_*=u$ almost everywhere in $\Omega_T$ furnishes the lower semi-continuity of $u$. In order to show this equality, we also define the set
\begin{equation}\label{Lebesgue}
\L= \bigg{\{} (x,t) \in \Omega_T : \, \, |u(x,t)|< \infty, \, \text{and} \, \lim_{\rho \rightarrow 0^+} \dashint_{(x,t)+\Q_{\rho}} |u(x,t)-u(y,t)|\, dydt=0\bigg{\}}.
\end{equation}
This set is well defined, since $u \in L^1_{loc}(0,T; L^{1}_{loc}(\Omega))$. Moreover,
\begin{equation}
\label{full}
|\L|=|\Omega_T|.
\end{equation}
As we will see, this is a consequence of the fact that $X:=(\Omega_T,\l^{N+1},d)$, being $\l^{N+1}$ the $(N+1)$-Lebesgue measure and $d$ a particular distance to be introduced, is a doubling space. We consider the following distance $d$: for any $(x,t), (y,s)\in\Omega_T$ we define 
\[d((x,t), (y,s)):= \max \{|x_i-y_i|^{\frac{p_i}{\bar{p}}}, |t-s|^{\frac{1}{\bar{p}}}\},\]
and we denote by $\b_{\rho}(x,t)$ the balls with respect to distance $d$. It turns out that $\b_{\rho}(x,t)=(x,t)+\Q_{\rho}$. The doubling property follows from
\[\l^{N+1}(\b_{2\rho})= \l^{N+1}(\Q_{2\rho})= (2\rho)^{N+\bar{p}}= 2^{N+\bar{p}}\l^{N+1}(\Q_{\rho})=2^{N+\bar{p}}\l^{N+1}(\b_{\rho}).\]
Accordingly, \cite[p. 12]{Heinonen} provides \eqref{full}. 
\newline
Taking \eqref{full} into account, it is sufficient to prove $u_*=u$ in $\L$. 
For all $(x,t) \in \L$ we have
\[u_*(x,t)=\lim_{\rho \rightarrow 0} \essinf_{(x,t)+\Q_{\rho}} u \leq \lim_{\rho \rightarrow 0} \dashint_{(x,t)+\Q_{\rho}} u \, dydt =u(x,t). \]
To show the opposite inequality, let us pick $(x_o,t_o) \in \L$ and suppose by contradiction that $u_*(x_o,t_o) < u(x_o,t_o)$. Let $r,b>0$ be small enough such that  $(x_o,t_o)+ \Q_r\subset \Omega_T$ and 
\[\essinf_{(x_o,t_o)+\Q_r} u := \mu_- \leq u_*(x_o,t_o)< \mu_-+ b<u(x_o,t_o).\] This choice is possible, since $u_*$ is close to $\mu^-$ for small values of $\rho$, as well as $\Q_\rho$ shrinks to $(x_o,t_o)$ as $\rho\to0^+$. Let us introduce $a \in (0,1)$ such that 
\[\mu_-+ab >u_*(x_o,t_o), \quad\quad \text{i.e.,} \quad\quad \frac{u_*(x_o,t_o)-\mu_-}{b} < a < 1. \]
\noindent Then there exists $\nu^-_a>0$, depending only on $a,b,p_i, N$, such that for some $\rho\in (0,r)$ we have
\[
|[u\leq\mu_-+b] \cap (x_o,t_o)+\Q_{\rho}| \leq \nu^-_a |\Q_{\rho}|,
\]
since otherwise we have, for all $\rho \in (0,r)$,
\[\begin{split}
\int_{(x_o,t_o)+\Q_{\rho}} |u(x_o,t_o)- u(x,t)| \, dxdt &\ge \int_{[u\leq\mu_-+b] \cap (x_o,t_o)+\Q_\rho} [u(x_o,t_o)-(\mu_-+b)]\, dxdt \\
&\ge \nu^-_a [u(x_o,t_o)-(\mu_-+b)] |\Q_{\rho}|,
\end{split}\]
contradicting $(x_o,t_o) \in \L$. Now we are in the position to apply Lemma \ref{DG} and reach
\[u(x,t) \ge \mu_-+ab>u_*(x_o,t_o), \quad \text{for a.a.}\, \, (x,t) \in (x_o,t_o)+\Q_{\rho/2}.\]
This contradicts the definition of $u_*(x_o,t_o)$, since 
\[u_*(x_o,t_o) <\essinf_{(x_o,t_o)+Q_{\rho/2}} u \leq \lim_{\rho \rightarrow 0^+} \essinf_{(x_o,t_o)+Q_{\rho}}u=u_*(x_o,t_o).\]
Therefore, we obtain $u_*=u$ in $\L$, concluding the proof.
\end{proof}

\begin{remark}
Semi-continuity of solutions to \eqref{EQ}, proved in Theorem \ref{semicontinuity}, jointly with the structure of the equation, ensures that the estimates \eqref{DG-TH} and \eqref{DG+TH} of Lemma \ref{DG} hold true everywhere (and not merely `almost everywhere') in $\Q_\rho^-(\theta)$. \end{remark}

% \section*{Appendix}\label{panini}
% \textcolor{blue}{Se uno non si fida del passaggio blu puo fare così.} Fisso $M= ||u||_{\infty, \Omega_T}$ all'inizio, lavoro solo su $\L_M(\Omega_T)$ e prendo $b \in (0,1)$ tale che 
% \[\inf_{Q_R(M)} u =\mu_-\leq u_*(0) < \mu_-+ bM< u(0).\] Adesso considero $a \in (0,1)$ tale che 
% \[ \mu_-+ abM> u_*(0) \quad \iff \quad \frac{u_*(0)-\mu_-}{bM}< a<1.\]

% Allora applico il De Giorgi Lemma \ref{DG} su 
% \[|[u<\mu_-+ bM] \cap \Q_R(M)| \leq \nu_-^a |\Q_R(M)|,\]
% per ottenere l'assurdo 
% \[u(0)> \mu_-+abM>u_*(0), \qquad \text{in} \quad Q_{R/2}(M).\]
% L'assurdo sta nel fatto che 
% \[u_*(0) =\lim_{\rho \rightarrow 0} \inf_{Q_{\rho}}u\]
% ma quando $\rho$ e' cosi piccolo che $Q_{\rho}$ entra dentro $Q_{R/2}(M)$ allora 
% \[u_*(0) <\inf_{Q_{R/2}(M)} u \leq \lim_{\rho \rightarrow 0} \inf_{Q_{\rho}}u=u_*(0).\]

% Nota che fare il minimo su un insieme piu piccolo e' piu' grande.
\section*{Acknowledgements}
The authors are grateful to Sunra J.N. Mosconi, for interesting conversations on the subject. S. Ciani acknowledges the support of the Department of Mathematics of T.U. Darmstadt. U. Guarnotta was supported by the following research projects: 1) PRIN 2017 ‘Nonlinear
Differential Problems via Variational, Topological and Set-valued Methods’ (Grant No. 2017AYM8XW) of
MIUR; 2) PRA 2020–2022 Linea 3 ‘MO.S.A.I.C.’ of the University of Catania. All the authors acknowledge the support of GNAMPA (INdAM).

%\newpage

%\section{Contazzi}
%Per la scelta di $\bar{n}$ in Lemma \ref{DG} basta richiedere che 
%\[\bigg(1+\frac{1}{2^{\bar{n}}} \bigg) \leq 2^{\bar{p}/p_N}, \quad \Rightarrow \quad 2^{\bar{n}} \ge \frac{1}{2^{\bar{p}/p_N}-1}, \quad \Rightarrow \quad \quad  \bar{n} \ge \log_2(1)-\log_2 (2^{\bar{p}/p_N}-1)=|\log_2(2^{\bar{p}/p_N}-1)|. \]
%Nota che sul tempo non c'e' problema perche' $\rho_0=2\rho$.\vskip0.1cm \noindent Per le cut-off vale
%\[\frac{1}{\rho_{n,i}-\rho_{n+1,i}}=\frac{1}{ \theta\rho (\frac{1}{2^{\bar{n}+n}})- \theta\rho (\frac{1}{2^{\bar{n}+n+1}})}=\frac{2^{\bar{n}+n+1}}{\theta \rho}=\gamma(\bar{n})\frac{2^{n}}{\theta \rho} .\]
%Mentre per i volumi vale
%\[|Q_0|= \bigg[\prod_{i=1}^N \theta^{f(i)} \rho^{g(i)} \bigg(1+\frac{1}{2^{\bar{n}}} \bigg)\bigg] (2\rho^{\bar{p}})=\bigg[\prod_i \theta^{f(i)} (2\rho)^{g(i)} \bigg] (2\rho^{\bar{p}}) 2^{-\sum_i g(i)} \bigg(1+\frac{1}{2^{\bar{n}}} \bigg)^N=\]
%\[ =|Q_{2\rho}^-(\theta)| \bigg(1+\frac{1}{2^{\bar{n}}} \bigg)^N2^{-\sum_i g(i)} =\gamma(\bar{n}) \, |Q_{2\rho}^-(\theta)|.\]
%\[\quad\text{con} \quad  {-\sum_i g(i)} \, \, \text{non dip. da i}\]

\end{document}